\input amstex
\input amsppt.sty
\magnification=\magstep1
\vsize=22.2truecm
\baselineskip=16truept
\pageno=1
\nologo
\TagsOnRight
\def\Z{\Bbb Z}

\def\Q{\Bbb Q}

\def\l{\left}
\def\r{\right}
\def\bg{\bigg}
\def\({\bg(}
\def\[{\bg\lfloor}
\def\){\bg)}
\def\]{\bg\rfloor}
\def\t{\text}
\def\f{\frac}

\def\p{\ (\roman{mod}\ p)}

\def\se {\subseteq}

\def\bi{\binom}
\def\eq{\equiv}

\def\ls{\leqslant}

\def\mo{\roman{mod}}

\def\ve{\varepsilon}
\def\al{\alpha}

\def\m#1#2{\thickfracwithdelims\{\}\thickness0{#1}{#2}_m}

\def\Proof{\noindent{\it Proof}}

\def\Remark{\medskip\noindent{\it  Remark}}
\def\Ack{\medskip\noindent {\bf Acknowledgment}}
\hbox {Accepted by Proc. Amer. Math. Soc. on March 6, 2005.}
\medskip
\topmatter
\title Binomial Coefficients and Quadratic Fields\endtitle
\author Zhi-Wei Sun\endauthor
\affil Department of Mathematics, Nanjing University\\
Nanjing 210093, People's Republic of China
\\zwsun\@nju.edu.cn
\ \ {\tt http://pweb.nju.edu.cn/zwsun}
\endaffil
\abstract Let $E$ be a real quadratic field with discriminant
$d\not\eq0\ (\mo\ p)$ where $p$ is an odd prime.
For $\rho=\pm1$ we determine
$\prod_{0<c<d,\ (\f dc)=\rho}\bi{p-1} {\lfloor pc/d\rfloor}$
modulo $p^2$ in terms of Lucas numbers, the
fundamental unit and the class number of $E$.
\endabstract
\thanks 2000 {\it Mathematics Subject Classifications}. Primary 11B65;
Secondary 11B37, 11B68, 11R11.
\newline\indent
The author was supported by the National Science Fund for
Distinguished Young Scholars (No. 10425103) and the Key Program of
NSF (No. 10331020) in China.
\endthanks
\endtopmatter
\document

\heading{1. Introduction}\endheading

Let $p$ be an odd prime not dividing a positive
integer $m$. A. Granville [G, (1.15)] discovered the remarkable congruence
$$\prod_{0<k<m}\bi{p-1}{\lfloor
pk/m\rfloor}\eq(-1)^{(m-1)(p-1)/2}(m^p-m+1)\ \ (\mo\ p^2),$$
where we use $\lfloor x\rfloor$ to denote the integral part of a real number $x$.
Subsequently the present author [S1] determined further
$\prod_{0<k<m/2}\bi{p-1}{\lfloor pk/m\rfloor}\ \mo\ p^2$.
In this paper a more sophisticated result connected
with real quadratic fields will be established.

 For $A,B\in\Z$ the Lucas sequences $u_n=u_n(A,B)$ and $v_n=v_n(A,B)$
($n=0,1,2,\ldots)$ are given by
$$\align&u_0=0,\ u_1=1,\ \t{and}\ u_{n+1}=Au_n-Bu_{n-1}
\ \t{for}\ n=1,2,3,\ldots,
\\&v_0=2,\ v_1=A,\ \t{and}\ v_{n+1}=Av_n-Bv_{n-1}\ \, \t{for}\ n=1,2,3,\ldots.
\endalign$$
It is well known that
$$(\al-\beta)u_n=\al^n-\beta^n\ \ \t{and}\ \
v_n=\al^n+\beta^n\quad\t{for every}\ n=0,1,2,\ldots,$$
where $\al$ and $\beta$ are the two roots of the equation $x^2-Ax+B=0$.
Also, for any odd prime $p$ we have
$u_p\eq(\f{\Delta}p)\ (\mo\ p)$ and $v_p\eq A\ (\mo\ p)$,
where $\Delta=A^2-4B$ and $(\f {\cdot}p)$
denotes the Legendre symbol. (See, e.g., [R, pp.\,41-55].)
If $p$ is an odd prime not dividing $B$, then
$p\mid u_{p-(\f{\Delta}p)}$ since $Au_p+v_p=2u_{p+1}$ and $Au_p-v_p=2Bu_{p-1}$.

 Throughout this paper, for an assertion $P$ we set
$$[P]=\cases1&\t{if}\ P\ \t{holds},\\0&\t{otherwise}.
\endcases\tag1.1$$

Our main result is as follows.

\proclaim{Theorem 1.1} Let $E$ be a quadratic field with discriminant
$d=2^{\al}p_1\cdots p_r$ where
$\al\in\{0,2,3\}$ and $p_1,\ldots,p_r$ are distinct odd primes.
Let $\ve=(a+b\sqrt d)/2$ be the fundamental unit of the field $E$
where $a,b\in\Z$, and $N(\ve)$ be the norm $(a^2-b^2d)/4$ of
$\ve$ with respect to the field extension $E/\Q$.
Let $h$ be the class number of the field $E$, and $p$ be an odd
prime not dividing $d$.
Then,
for $\rho =\pm1$ we have
$$\aligned\prod\Sb 0<c<d\\(\f dc)=\rho  \endSb\bi{p-1}{\lfloor
pc/d\rfloor}\eq& 1+\f{\varphi(d)}2\((\al+[\al>0])(2^{p-1}-1)
+\sum_{0<i\ls r}\f{p_i^p-p_i}{p_i-1}\)
\\&+\f{\rho }2\l(\f dp\r)^{[N(\ve)=1]}
u_{p-(\f dp)}(a,N(\ve))bdh\ \ (\mo\ p^2),
\endaligned\tag1.2$$
where $\varphi$ is Euler's totient function and
$(\f d{\cdot})$ is the Kronecker symbol.
\endproclaim

\Remark. Under the conditions of Theorem 1.1,
$d\eq1\ (\mo\ 4)$ if $\al=0$, and $d/4\eq 3\ (\mo\ 4)$ if $\al=2$;
also $p$ divides
$bu_{p-(\f dp)}(a,N(\ve))$ since for $p\nmid b$ we have
$$\(\f{a^2-4N(\ve)}p\)=\(\f{b^2d}p\)=\l(\f dp\r).$$

\medskip
\noindent{\it Example}.
Each of the quadratic fields $\Q(\sqrt{13}),\Q(\sqrt{21}),
\Q(\sqrt6),\Q(\sqrt7)$
has class number 1, and their fundamental units are
$$\f{3+\sqrt{13}}2,\ \f{5+\sqrt{21}}2,\
5+2\sqrt6=\f{10+2\sqrt{24}}2,\ 8+3\sqrt7=\f{16+3\sqrt{28}}2$$
with norms $-1,1,1,1$ respectively; see, e.g., [C, p.\,271].
Let $p$ be an odd prime and $\rho \in\{1,-1\}$.
If $p$ does not divide 13, 21, 6, and 7, respectively, then
Theorem 1.1 gives the congruences
$$\align\prod\Sb 0<c<13\\(\f{13}c)=\rho  \endSb\bi{p-1}{\lfloor
pc/13\rfloor}\eq&1+\f{13^p-13}2+\rho  \f{13}2u_{p-(\f{13}p)}(3,-1),
\\\prod\Sb 0<c<21\\(\f{21}c)=\rho  \endSb\bi{p-1}{\lfloor
pc/21\rfloor}\eq&1+3(3^p-3)+7^p-7+
\rho\l(\f{21}p\r)\f{21}2u_{p-(\f{21}p)}(5,1),
\\\prod\Sb 0<c<24\\2\nmid c,\ (\f{6}c)=\rho  \endSb\bi{p-1}{\lfloor
pc/24\rfloor}\eq&1+8(2^p-2)+2(3^p-3)+\rho\l(\f 6p\r)24u_{p-(\f{6}p)}(10,1),
\\\prod\Sb 0<c<28\\2\nmid c,\ (\f{7}c)=\rho  \endSb\bi{p-1}{\lfloor
pc/28\rfloor}\eq&1+9(2^p-2)+7^p-7+\rho  \l(\f 7p\r)42u_{p-(\f{7}p)}(16,1)
\endalign$$
modulo $p^2$ respectively, where $(\f 6c)$ and $(\f 7c)$ are Jacobi symbols.
\medskip

We deduce Theorem 1.1 by combining the following two theorems.

\proclaim{Theorem 1.2} Let $m>2$ be an integer with the
factorization $p_1^{\al_1}\cdots p_r^{\al_r}$ where $p_1,\ldots,p_r$
are distinct primes and $\al_1,\ldots,\al_r$ are positive integers.
Let $p$ be an odd prime not dividing $m$. Then
$$\aligned&(-1)^{\f{\varphi(m)}2\cdot\f{p-1}2}\l(\f{p_1}p\r)^{[r=1]}
\prod\Sb 0<k<m/2\\(k,m)=1\endSb\bi{p-1}{\lfloor pk/m\rfloor}
\\\eq&1+\f{\varphi(m)}2\sum_{i=1}^r(\al_ip_i-\al_i+1)\f{p_i^{p-1}-1}{p_i-1}
\ \ (\mo\ p^2).
\endaligned\tag1.3$$
\endproclaim

 In the next theorem we use the Bernoulli polynomial $B_n(x)$ of degree $n$
 and the $n$th Bernoulli number $B_n=B_n(0)$.
 Also, we let $\Bbb P$ denote the set of all (positive) primes.

\proclaim{Theorem 1.3} Let $E$ be a real quadratic field with
discriminant $d$ and class number $h$. Let
$\ve=(a+b\sqrt d)/2>1$ be the fundamental unit of $E$ where $a,b\in\Z$,
and $N(\ve)$ be the norm $(a^2-b^2d)/4$ of $\ve$. Let $p$
be an odd prime not dividing $d$,
and let $u$ stand for $bu_{p-(\f dp)}(a,N(\ve))$.
Then
$$\sum_{c=1}^{d-1}\l(\f dc\r)\l(B_{p-1}\l(\f cd\r)-B_{p-1}\r)
\eq\l(\f dp\r)^{[N(\ve)=-1]}dh\f{u}p\ \  (\mo\ p),\tag1.4$$
and
$$\prod\Sb 0<c<d/2\\(c,d)=1\endSb\bi{p-1}{\lfloor
pc/d\rfloor}^{(\f dc)}
\eq\cases(\f dp)(1+\f{dhu}2)\ (\mo\ p^2)&\t{if}\ d=8\ \t{or}\ d\in\Bbb P,
\\1+(\f dp)^{[N(\ve)=1]}\f{dhu}2\ (\mo\ p^2)&\t{otherwise}.\endcases\tag1.5$$
\endproclaim

\Remark. In the case where $d\eq 1\ (\mo\ 4)$ is a prime,
(1.4) was proved in [GS] by means of $p$-adic
logarithms and  Dirichlet's class number formula (see, e.g., [W]).
\medskip

 In the spirit of R. Crandall and C. Pomerance [CP],
Theorems 1.1--1.3 might be of computational interest.

We shall make some preparations in the next section
and give proofs of Theorems 1.1--1.3
in Section 3.

\heading{2. On the sum $\sum\Sb 0<k<p\\m\mid k-r\endSb\f1k$ modulo $p$}
\endheading

Bernoulli polynomials play important roles in many
aspects. The reader is referred to [IR, pp.\,228-248] for basic properties,
and to [DSS] for a bibliography of related papers.

In this section we prove the following basic result
and derive some consequences.

\proclaim{Theorem 2.1} Let $m$ be a positive integer
not divisible by an odd prime $p$. Then
for any $r\in\Z$ we have
$$\sum^{p-1}\Sb k=1\\k\eq r\ (\mo\ m)\endSb\f1k
\eq\f1m\(B_{p-1}\l(\l\{\f rm\r\}\r)-B_{p-1}\l(\l\{\f{r-p}m\r\}\r)\)
\ \ (\mo\ p),\tag2.1$$
where $\{x\}$ stands for the fractional part of a real number $x$.
\endproclaim
\Proof. Applying Lemma 3.1 of [S3] with $k=p-2$, we find that
$$-m\sum^{p-1}\Sb j=1\\j\eq r\ (\mo\ m)\endSb\f1j\eq B_{p-1}\l(\f pm
+\l\{\f{r-p}m\r\}\r)-B_{p-1}\l(\l\{\f{r}m\r\}\r)\ (\mo\ p).$$
For $t=\{(r-p)/m\}$, we have
$$B_{p-1}\l(\f pm+t\r)-B_{p-1}(t)
=\sum_{l=1}^{p-1}\bi{p-1}lB_{p-1-l}\l(\l(\f pm+t\r)^l-t^l\r)
\eq0\ (\mo\ p).$$
(Recall that $B_1=-1/2$ and $B_{2n+1}=0$ for $n=1,2,\ldots$.
Also, $p$ divides no denominators of $B_0, B_2,\ldots,B_{p-3}$
by the theorem of Clausen and von Staudt (cf. [IR, pp.\,233-236]).)
Therefore (2.1) follows.
\qed

\Remark. The author first discovered Theorem 2.1 in Sept. 1991
by using Fourier series, and Lemma 3.1 of [S3] was originally motivated by
this result.

\proclaim{Corollary 2.1} Let $m$ and $n$ be positive integers,
and let $p$ be an odd prime not dividing $m$. Then
$$\quad\ B_{p-1}\l(\l\{\f{pn}m\r\}\r)-B_{p-1}
\eq m\sum_{r=1}^n K_p(r,m)
\eq-\sum^{\lfloor pn/m\rfloor}\Sb k=1\\p\nmid k\endSb\f1k
\ \ (\mo\ p),
\tag2.2$$
where
$$K_p(r,m):=\sum^{p-1}\Sb k=1\\m\mid k-rp\endSb\f1k
=\sum^{p-1}\Sb l=1\\m\mid l-(1-r)p\endSb\f1{p-l}
\eq-K_p(1-r,m)\ (\mo\ p).\tag2.3$$
\endproclaim

\Proof. In view of Theorem 2.1,
$$\align m\sum_{r=1}^n K_p(r,m)\eq&\sum_{r=1}^n\(B_{p-1}\l(\l\{\f{rp}m\r\}\r)
-B_{p-1}\l(\l\{\f{(r-1)p}m\r\}\r)\)
\\\eq&B_{p-1}\l(\l\{\f{pn}m\r\}\r)-B_{p-1}\ (\mo\ p).
\endalign$$
On the other hand,
$$-\sum_{r=1}^nK_p(r,m)\eq \sum_{r=1}^n\sum^{p-1}
\Sb k=1\\ m\mid rp-k\endSb\f{1}{rp-k}
=\sum^{pn-1}\Sb j=1\\ p\nmid j,\ m\mid j\endSb\f{1}j
=\sum^{\lfloor pn/m\rfloor}\Sb k=1\\p\nmid k\endSb\f{1}{km}
 \ (\mo\ p).$$
So we have (2.2).\qed

Let $p$ be an odd prime and $r$ be any integer.
An explicit congruence for $K_p(r,12)$ mod $p$
appeared in Corollary 3.3 of [S2].
By Theorem 2.1 and [GS, (4)]
we can also determine
$$K_p(3+6r,24),\ K_p(5,40),\ K_p(25,40),\ K_p(6,60),\ K_p(36,60)$$
modulo $p$ in terms of some second-order linear recurrences.

 For a prime $p$ and any $a\in\Z$ not divisible by $p$,
the Fermat quotient $q_p(a)$ is defined as the integer $(a^{p-1}-1)/p$.

\proclaim{Corollary 2.2} Let $p$ be an odd prime
and let $m$ be a positive integer
not divisible by $p$. Then we have
$$\sum_{r=1}^mrK_p(r,m)\eq-q_p(m)\ \ (\mo\ p).\tag 2.4$$
\endproclaim
\Proof. By Corollary 2.1,
$$\align\sum_{n=1}^m\sum_{r=1}^nK_p(r,m)\eq&\f1 m\sum_{n=1}^m
\(B_{p-1}\l(\l\{\f{pn}m\r\}\r)-B_{p-1}\)
\\\eq&m^{p-2}\(\sum_{n=1}^mB_{p-1}\l(\l\{\f {pn}m\r\}\r)-mB_{p-1}\)
\\\eq&\sum_{r=0}^{m-1}m^{p-2}B_{p-1}\l(\f rm\r)-m^{p-1}B_{p-1}
=(1-m^{p-1})B_{p-1}\ (\mo\ p)
\endalign$$
where we have applied Raabe's theorem in the last step.
It is well known that $pB_{p-1}\eq-1\ (\mo\ p)$ (cf. [IR, p.\,233]). Also,
$$\align&\sum_{n=1}^m\sum_{r=1}^nK_p(r,m)
=\sum_{r=1}^m(m-(r-1))K_p(r,m)
\\\eq&-\sum_{r=1}^m(m+1-r)K_p(m+1-r,m)=-\sum_{s=1}^msK_p(s,m)\ (\mo\ p).
\endalign$$
So we have (2.4). \qed

\Remark. It can be shown that (2.4) is equivalent to
a formula of Lerch [L] which was deduced in a different way.

\heading{3. Proofs of Theorems 1.1--1.3}\endheading

\noindent {\it Proof of Theorem 1.2}. For each positive integer $d$ we set
$$\psi(d)=\prod\Sb 0<c<d/2\\(c,d)=1\endSb\bi{p-1}{\lfloor pc/d\rfloor},$$
where $\psi(1)$ and $\psi(2)$ are considered as 1. For any $a\in\Z$
with $p\nmid a$, clearly
$$\align a^p-a=&a\l(a^{(p-1)/2}+\l(\f ap\r)\r)\l(a^{(p-1)/2}-\l(\f ap\r)\r)
\\\eq &2a\l(\f ap\r)\l(a^{(p-1)/2}-\l(\f ap\r)\r)\ (\mo\ p^2).
\endalign$$
Thus, Theorem 1.1 of [S1] implies that if $d\not\eq0\ (\mo\ p)$ then
$$\align &(-1)^{\f{p-1}2\l\lfloor\f{d-1}2\r\rfloor}
\prod_{0<c<d/2}\bi{p-1}{\lfloor pc/d\rfloor}
\\\eq&\cases(\f dp)+(\f dp)\f{d^p-d}2&\t{if}\ 2\nmid d,
\\(\f{2d}p)+(\f{2d}p)\f{d^p-d}2-(\f{2d}p)\f{2^p-2}2&\t{if}\ 2\mid d,
\endcases
\\\eq&\l(\f dp\r)\l(\f 2p\r)^{[2\mid d]}
\l(1+\f{d^p-d}2-[2\mid d](2^{p-1}-1)\r)\ \ (\mo\ p^2).
\endalign$$
Since $\prod_{0<k<n/2}\bi{p-1}{\lfloor pk/n\rfloor}=\prod_{d\mid n}\psi(d)$
for $n=1,2,\ldots$,
applying the M\"obius inversion formula we get that
$$\align &\psi(m)=\prod_{d\mid m}\prod_{0<c<d/2}
\bi{p-1}{\lfloor pc/d\rfloor}^{\mu(m/d)}
\\\eq&(-1)^{\f{p-1}2\sum_{d\mid m}\mu(\f md)(\f{d-1}2-\f{[2\mid d]}2)}
\l(\f 2p\r)^{\sum_{d\mid m}\mu(m/d)[2\mid d]}
\\&\times\prod_{d\mid m}\l(\f dp\r)^{\mu(m/d)}
\times\prod_{d\mid m}\l(1+\mu\l(\f md\r)\l(\f{d^p-d}2-[2\mid d](2^{p-1}-1)\r)\r)
\  (\mo\ p^2).
\endalign$$
By elementary number theory, $\sum_{d\mid m}\mu(\f md)\f{d-1}2=\f{\varphi(m)}2$
and also
$$\sum_{d\mid m}\mu\l(\f md\r)[2\mid d]=\sum_{2c\mid m}\mu\l(\f{m}{2c}\r)
=[2\mid m]\sum_{c\mid (m/2)}\mu\l(\f{m/2}c\r)=0$$
since $m>2$. Therefore
$$(-1)^{\f{\varphi(m)}2\cdot\f{p-1}2}\psi(m)
\eq\prod_{d\mid m}\l(\f dp\r)^{\mu(m/d)}
\times\bigg(1+\sum_{d\mid m}\mu\l(\f md\r)\f{d^p-d}2\bigg)
\ (\mo\ p^2).$$

Observe that
$$\align&\prod_{d\mid m}\l(\f
dp\r)^{\mu(m/d)}=\prod_{I\se\{1,\ldots,r\}}
\l(\f{m/\prod_{i\in I}p_i}p\r)^{\mu(\prod_{i\in I}p_i)}
\\=&\l(\f{m^{2^r}/\prod_{I\se\{1,\ldots,r\}}\prod_{i\in I}p_i}p\r)
=\l(\f{m^{2^r}/\prod_{i=1}^rp_i^{2^{r-1}}}p\r)
\\=&\l(\f{\prod_{i=1}^rp_i^{2^{r-1}(2\al_i-1)}}p\r)
=\l(\f{p_1\cdots p_r}p\r)^{2^{r-1}}=\l(\f{p_1}p\r)^{[r=1]}.
\endalign$$
Also,
$$\align&\varphi(m)+\sum_{d\mid m}\mu\l(\f md\r)(d^p-d)
=\sum_{d\mid m}\mu(d)\f{m^p}{d^p}
=m^p\prod_{i=1}^r\l(1-p_i^{-p}\r)
\\=&\prod_{i=1}^r\l(p_i^{\al_i p}-p_i^{(\al_i-1)p}\r)
=\prod_{i=1}^r\l((p_i+(p_i^p-p_i))^{\al_i}-(p_i+(p_i^p-p_i))^{\al_i-1}\r)
\\\eq&\prod_{i=1}^r\l(p_i^{\al_i}+\al_i
p_i^{\al_i-1}(p_i^p-p_i)-\l(p_i^{\al_i-1}
+(\al_i-1)p_i^{\al_i-2}(p_i^p-p_i)\r)\r)
\\\eq&\prod_{i=1}^r\l(\varphi(p_i^{\al_i})+(p_i^{p-1}-1)
(\al_ip_i^{\al_i}-(\al_i-1)p_i^{\al_i-1})\r)
\\\eq&\varphi(m)\(1+\sum_{i=1}^r\f{p_i^{p-1}-1}{p_i-1}
(\al_i p_i-\al_i+1)\)\ (\mo\ p^2).
\endalign$$
Thus (1.3) holds in view of the above. \qed

\medskip
\noindent{\it Proof of Theorem 1.3}.
Write $\ve^{p-(\f dp)}=(V+U\sqrt d)/2$ where $U,V\in\Z$,
and let $p'$ be an integer with $pp'\eq1\ (\mo\ d)$.
Theorem 3.1 of Williams [W] states that
$$h\f Up\eq-\l(\f dp\r)N(\ve)^{((\f dp)-1)/2}
\sum_{i=1}^{p-1}\f{\beta_p(i)}i\ \ (\mo\ p)$$
where $\beta_p(i)=\sum_{0<j<d\{p'i/d\}}(\f dj)$.

Let $\bar\ve=(a-b\sqrt d)/2$. Then $\ve+\bar\ve=a$ and
$\ve\bar\ve=N(\ve)$. Clearly
$$v_n(a,N(\ve))+u_n(a,N(\ve))b\sqrt d
=\ve^n+{\bar\ve}^n+\f{\ve^n-{\bar\ve}^n}{\ve-\bar\ve}b\sqrt d=2\ve^n$$
for $n=0,1,\ldots$, thus
$U=bu_{p-(\f dp)}(a,N(\ve))=u$ (and $V=v_{p-(\f dp)}(a,N(\ve))$).

Observe that
$$\align\sum_{i=1}^{p-1}\f{\beta_p(i)}i=
&\sum_{j=1}^{d-1}\l(\f dj\r)\sum\Sb0<i<p\\d\{p'i/d\}>j\endSb\f1i
=\sum_{j=1}^{d-1}\l(\f dj\r)\sum_{j<r<d}
\sum\Sb 0<i<p\\d\mid p'i-r\endSb\f1i
\\=&\sum_{j=1}^{d-1}\l(\f dj\r)\sum_{j<r<d}
\sum\Sb 0<i<p\\d\mid i-rp\endSb\f1i
=\sum_{j=1}^{d-1}\l(\f dj\r)\sum_{j<r<d}K_p(r,d).
\endalign$$
As $\chi(j)=(\f dj)$ is a nontrivial multiplicative character
modulo $d$, the sum $\sum_{j=1}^{d-1}(\f dj)$ vanishes.
Therefore, with the help of Corollary 2.1, we have
$$\align \sum_{i=1}^{p-1}\f{\beta_p(i)}i
=&\sum_{j=1}^{d-1}\l(\f dj\r)\(\sum_{r=1}^dK_p(r,d)-\sum_{r=1}^jK_p(r,d)\)
\\\eq&\sum_{j=1}^{d-1}\l(\f dj\r)
\f1d\l(0-B_{p-1}\l(\l\{\f{pj}d\r\}\r)+B_{p-1}\r)
\\\eq&-\f1d\l(\f dp\r)\sum_{j=1}^{d-1}
\l(\f d{pj}\r)\l(B_{p-1}\l(\l\{\f{pj}d\r\}\r)-B_{p-1}\r)
\\\eq&-\f1d\l(\f dp\r)\sum_{c=1}^{d-1}\l(\f dc\r)
\l(B_{p-1}\l(\f cd\r)-B_{p-1}\r)
\ \ (\mo\ p).\endalign$$

Combining the above we obtain (1.4).

For each $c=1,\ldots,d-1$, we have
$\chi(d-c)=\chi(-1)\chi(c)=\chi(c)$; also
$$\align &(-1)^{\lfloor pc/d\rfloor}\bi{p-1}{\lfloor pc/d\rfloor}
=\prod_{k=1}^{\lfloor pc/d\rfloor}\l(1-\f pk\r)
\\\eq&1-p\sum_{k=1}^{\lfloor pc/d\rfloor}\f1k\eq
1+p\(B_{p-1}\l(\l\{\f{pc}d\r\}\r)-B_{p-1}\)\ (\mo\ p^2).
\endalign$$
Taking the above congruence and (1.3) modulo $p$, we obtain
$$\align\prod\Sb 0<c<d/2\\(c,d)=1\endSb(-1)^{\lfloor pc/d\rfloor}
\eq&\prod\Sb 0<c<d/2\\(c,d)=1\endSb\bi{p-1}{\lfloor pc/d\rfloor}
\\\eq&(-1)^{\f{\varphi(d)}2\cdot\f{p-1}2}\l(\f dp\r)^{[d\
\t{is a prime power}]}\ \ (\mo\ p)
\endalign$$
and hence
$$\prod_{0<c<d/2}(-1)^{\lfloor pc/d\rfloor(\f dc)}=\l(\f
dp\r)^{[d=8\ \t{or}\ d\in\Bbb P]}.$$
(Note that $4\mid\varphi(d)$ and no square of
an odd prime divides $d$.)
On the other hand,
$$\align&\prod_{0<c<d/2}\((-1)^{\lfloor pc/d\rfloor}
\bi{p-1}{\lfloor pc/d\rfloor}\)^{(\f dc)}
\\\eq&\prod_{0<c<d/2}\(1+p\l(\f dc\r)
\l(B_{p-1}\l(\l\{\f{pc}d\r\}\r)-B_{p-1}\r)\)
\\\eq&1+\f p2\sum_{0<c<d/2}\l(\f dc\r)
\l(B_{p-1}\l(\l\{\f{pc}d\r\}\r)-B_{p-1}\r)
\\&+\f p2\sum_{0<c<d/2}\l(\f d{d-c}\r)\l(B_{p-1}
\l(\l\{\f{p(d-c)}d\r\}\r)-B_{p-1}\r)
\\\eq&1+\f p2\sum_{c=1}^{d-1}\l(\f dc\r)
\l(B_{p-1}\l(\l\{\f{pc}d\r\}\r)-B_{p-1}\r)
\\\eq&1+\f p2\l(\f dp\r)\sum_{r=1}^{d-1}\l(\f dr\r)
\l(B_{p-1}\l(\f rd\r)-B_{p-1}\r)
\ \ (\mo\ p^2).
\endalign$$
These, together with (1.4), yield
$$\prod_{0<c<d/2}\bi{p-1}{\lfloor pc/d\rfloor}^{(\f dc)}
\eq\l(\f dp\r)^{[d=8\ \t{or}\ d\in\Bbb P]}\(1+\f{dhu}2\l(\f dp\r)^{[N(\ve)=1]}\)
\ \ (\mo\ p^2).$$
It is well known that $N(\ve)=-1$ if $d=8$ or $d\in\Bbb P$
(see, e.g., [C, pp.\,185-186]). So the desired (1.5) follows. \qed

\medskip
\noindent{\it Proof of Theorem 1.1}.
By Theorem 1.2 and the proof of Theorem 1.3,
$$\l(\f dp\r)^{[d=8\ \t{or}\ d\in\Bbb P]}\prod\Sb
0<c<d/2\\(c,d)=1\endSb\bi{p-1}{\lfloor pc/d\rfloor}
\eq1+\f{\varphi(d)}2F(d,p)
\ \ (\mo\ p^2)$$
where
$$\align F(d,p)=&[\al>0](2\al-\al+1)\f{2^{p-1}-1}{2-1}
+\sum_{0<i\ls r}(p_i-1+1)\f{p_i^{p-1}-1}{p_i-1}
\\=&(\al+[\al>0])(2^{p-1}-1)+\sum_{0<i\ls r}\f{p_i^p-p_i}{p_i-1};
\endalign$$
also
$$\l(\f dp\r)^{[d=8\ \t{or}\ d\in\Bbb P]}\prod\Sb
0<c<d/2\\(c,d)=1\endSb\bi{p-1}{\lfloor pc/d\rfloor}^{(\f dc)}
\eq1+\f{dhu}2\l(\f dp\r)^{[N(\ve)=1]}\ \ (\mo\ p^2)$$
where $u=bu_{p-(\f dp)}(a,N(\ve))\eq0\ (\mo\ p)$.
Therefore
$$\align&\prod\Sb 0<c<d/2\\(\f dc)=\rho \endSb\bi{p-1}{\lfloor pc/d\rfloor}
\bi{p-1}{\lfloor p(d-c)/d\rfloor}
=\prod\Sb 0<c<d/2\\(c,d)=1\endSb\bi{p-1}{\lfloor pc/d\rfloor}^{1+\rho(\f dc)}
\\\eq&\l(1+\f{\varphi(d)}2F(d,p)\r)\(1+\f{dhu}2\l(\f dp\r)^{[N(\ve)=1]}\)^{\rho}
\\\eq&\l(1+\f{\varphi(d)}2F(d,p)\r)\(1+\rho \f{dhu}2\l(\f dp\r)^{[N(\ve)=1]}\)
\\\eq&1+\f{\varphi(d)}2F(d,p)+\rho\f{dhu}2\l(\f dp\r)^{[N(\ve)=1]}
\ \ (\mo\ p^2).
\endalign$$
This proves (1.2). We are done. \qed

\Ack. The author thanks the referee for his many helpful comments.

\enddocument